\documentclass[12pt,dvips]{amsart}
\usepackage{euler, amsfonts, amssymb, latexsym, epsfig,epic}

\newcommand\B{\reals^3_{\sum=0}}

\newcommand\iso{{\,\cong\,}}
\newcommand\tensor{{\otimes}}

\newtheorem{Theorem}{Theorem} 
\newtheorem{Proposition}{Proposition} 
\newtheorem{Lemma}{Lemma}

\newtheorem*{Theorem*}{Theorem}
\newtheorem*{Lemma*}{Lemma}

\newtheorem*{Example*}{Example}
\newtheorem*{Proposition*}{Proposition}

\newcommand\union{\bigcup}

\newcommand\reals{{\mathbb R}}
\newcommand\complexes{{\mathbb C}}
\newcommand\integers{{\mathbb Z}}

\newcommand\onehalf{\frac{1}{2}}
\newcommand\naturals{{\mathbb N}}

\newcommand\dec{\hbox{\ \rm\tiny dec}}

\newcommand\Oplus{\bigoplus}

\newcommand\GLn{GL_n(\complexes)}

\theoremstyle{plain}

\newcommand\dfn{\bf} 

\begin{document}
\pagestyle{plain}

\title{A positive proof of the Littlewood-Richardson rule \\
  using the octahedron recurrence}
\author{Allen Knutson}
\email{allenk@math.berkeley.edu}
\thanks{AK was supported by NSF grant 0072667, and a Sloan Fellowship.}
\address{Mathematics Department\\ UC Berkeley\\ Berkeley, California}
\author{Terence Tao}
\thanks{TT was supported by the Clay Mathematics Institute,
  and the Packard Foundation.} 
\email{tao@math.ucla.edu}
\address{Mathematics Department\\ UCLA\\ Los Angeles, California}
\author{Christopher Woodward}
\thanks{CW was supported by NSF grant 9971357.}
\email{ctw@math.rutgers.edu}
\address{Mathematics Department\\ Rutgers University\\ New Brunswick, New Jersey}
\date{\today}

\maketitle

\begin{abstract}
  We define the {\em hive ring}, which has a basis indexed by dominant weights
  for $\GLn$, and structure constants given by counting {\em hives} 
  \cite{KT1} (or equivalently honeycombs, or Berenstein-Zelevinsky patterns
  \cite{BZ1}). We use the octahedron rule from \cite{RR,FZ,P,S} to prove 
  bijectively that this ``ring'' is indeed associative.

  This, and the Pieri rule, give a self-contained proof that 
  the hive ring is isomorphic as a ring-with-basis
  to the representation ring of $\GLn$.
  
  In the honeycomb interpretation, the octahedron rule becomes
  ``scattering'' of the honeycombs. This recovers some of the
  ``crosses and wrenches'' diagrams from the very recent preprint \cite{S},
  whose results we use to give a closed form for the associativity bijection.

\end{abstract}

\tableofcontents

\newcommand\Oone{{\mathcal O}(1)}
\newcommand\PP{{\mathbb P}}

\newcommand\HIVE{{\tt HIVE}}
\newcommand\HONEY{{\tt HONEY}}

\section{Introduction}

Let $Rep(\GLn)$ denote the ring of (formal differences of algebraic
finite-dimensional) representations of $\GLn$, with addition and
multiplication coming from direct sum and tensor product of representations.
Then $Rep(\GLn)$ has a canonical basis $\big\{[V_\lambda]\big\}$, 
the irreducible representations, 
indexed by the set $\integers^n_{\dec}$
of weakly decreasing $n$-tuples of integers.
(The ``$[\,\,]$'' are only there to maintain a proper distinction between
an actual representation $V_\lambda$ and the corresponding element of
$Rep(\GLn)$, which is really an isomorphism class.) Our reference for
this representation theory is \cite{FH}.

The structure constants $c_{\lambda\mu}^\nu$ of this ring-with-basis,
defined by 
$$ [V_\lambda] \, [V_\mu] = \sum_\nu c_{\lambda\mu}^\nu \, [V_\nu],$$
are necessarily nonnegative (being the dimensions of certain vector spaces
of intertwining operators), and there are many known rules for 
calculating them
as the cardinalities of certain
combinatorially defined sets. The most famous is the Littlewood-Richardson
rule, which counts skew Young tableaux.

In several of these rules, the set being counted is the lattice
points in a polytope (and in fact the polytopes from the
different rules are linearly equivalent). 
The first was in the unpublished thesis \cite{J}, and was proved by
establishing a bijection with skew Young tableaux; 
see also the appendix to \cite{B}.
It was rediscovered in \cite{BZ1}, where the proof starts with the
(nonpositive) Steinberg rule for tensor products and uses an
involution to cancel the negative terms. There is another extremely
roundabout proof via the connection with Schubert calculus, for which
a self-contained proof of a combinatorial rule was given in \cite{KT2}.

In this paper we give a new self-contained proof of this 
lattice-point-counting rule, 
in its incarnation as counting the {\em hives} from \cite{KT1},
whose definition we recall below. The main difficulty is in proving
that the ring so defined (which is supposed to match up with
$Rep(\GLn)$) is associative. We give a bijective proof of this,
using the octahedron rule from \cite{RR,P,FZ,S}.
This bijection was first found by CW in the honeycomb model, where
the connection to the octahedron rule is not transparent.

Very recently, in \cite{S}, a closed form was found for compositions
of the octahedron rule. In the last section we describe this formula
in the special case relevant for this paper. 

Since the octahedron rule is related to tropical algebraic geometry,
we hope that our bijective proof of associativity will turn out to be
the tropicalization of some natural but heretofore undiscovered
birational map, as in \cite{BZ3}.

\subsection{Acknowledgements.} 
It is our pleasure to thank Andrei Zelevinsky for comments
on an earlier version of this paper, David Speyer for kindly working out
the special case of his results \cite{S} which appears in section
\ref{sec:speyer}, and Jim Propp for directing us to Speyer's preprint.

\section{Hives}

Consider the triangle $\big\{ [x,y,z]: x$+$y$+$z=n$, $x,y,z\geq 0 \big\}$. 
This has $n+2 \choose 2$ integer points; call this finite set $tri_n$.
We will draw it in the plane and put $[n,0,0]$ at the lower left,
$[0,n,0]$ at the top, and $[0,0,n]$ in the lower right. 
This triangle breaks up into ${n+1 \choose 2}$ right-side-up triangles 
$\overline{[x+1,y,z] \,[x,y+1,z]\,[x,y,z+1]}$ and ${n\choose 2}$
upside-down triangles
$\overline{[x-1,y,z] \,[x,y-1,z]\,[x,y,z-1]}$.
We will count certain integer labelings of $tri_n$ to compute
Littlewood-Richardson coefficients, following \cite{J}, \cite{BZ1},
and especially \cite{KT1}.

\begin{figure}[htbp]
  \centering
  \epsfig{file=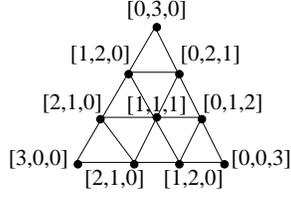,height=1in}
  \caption{The set $tri_3$, with its $3+1 \choose 2$ right-side-up and
    $3 \choose 2$ upside-down triangles.}
  \label{fig:tri_n}
\end{figure}

A {\dfn hive of size $n$} 
is a function $h: tri_n \to \integers$ satisfying certain inequalities. 
Here are three equivalent ways to state those inequalities (of which
we shall mainly use the first):
\begin{enumerate}
\item 
$ h_{x+1,y,z+1} + h_{x,y+1,z+1} \geq h_{x+1,y+1,z} + h_{x,y,z+2} $ 
when these four points are all in $tri_n$, and likewise for the
$120^\circ$ and $240^\circ$ rotations of the hive.
\item
If you extend $h$ to a real-valued function on the solid triangle
by making it linear on each little triangle 
$\overline{[x\pm 1,y,z] \,[x,y\pm 1,z]\,[x,y,z\pm 1]}$, $h$ is convex.
\item 
On each unit rhombus in the triangle, the sum across the short diagonal
is greater than or equal to the sum across the long diagonal.
\end{enumerate}
Note that the definition also makes sense for real-valued functions, 
in which case we will speak of a {\dfn real hive}. (We won't use
this concept until section \ref{sec:assoc}.)

Call these inequalities the {\dfn rhombus inequalities} on a hive.
They naturally come in three families, according to the orientation of
the rhombus.

\begin{figure}[htbp]
  \centering
  \epsfig{file=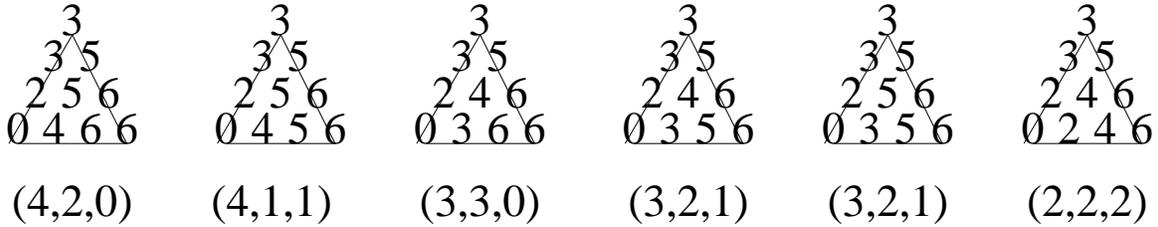,width=6in}
  \caption{The hives with Northwest and Northeast side 
    having differences $(2,1,0)$. The differences across the South side
    are indicated.}
  \label{fig:minisamplehives}
\end{figure}

Definitions linearly equivalent to this one appeared first in \cite{J,BZ1,BZ2}.
This version from \cite{KT1}, like the one in \cite{BZ2}, 
has the benefit that each inequality only involves a constant 
number of entries (namely four), independent of $n$.

\begin{Proposition}\label{prop:dec}
  Let $a_0,a_1,\ldots,a_n$ be the numbers on one side of a hive
  (read left-to-right).
  Then $a$ is convex, i.e. $a_i \geq \onehalf (a_{i-1}+a_{i+1})$.
  Put another way, the list $(a_1 - a_0, a_2-a_1,\ldots, a_n-a_{n-1})$ 
  is a weakly decreasing list of integers.
\end{Proposition}

\begin{proof}
  There are two rhombi with an obtuse vertex at $a_i$. Adding the two
  corresponding rhombus inequalities, we get the desired result.
\end{proof}

We can interpret such a list as a dominant weight for $\GLn$;
call the set of such weights $\integers^n_{\dec}$, 
and let $\lambda,\mu,\nu \in \integers^n_{\dec}$
be three of them.
Let $\HIVE_{\lambda\mu}^\nu$ denote the set of hives of size $n$ such that
\begin{itemize}
\item the lower left entry is zero
\item the differences on the Northwest side of the hive give $\lambda$
\item the differences on the Northeast side of the hive give $\mu$
\item the differences on the South side of the hive give $\nu$
\end{itemize}
where all differences are computed {\em left-to-right} throughout the paper. 
Note that for $\HIVE_{\lambda\mu}^\nu$ to be nonempty, we must have
$\sum_i (\lambda_i+\mu_i) = \sum_i \nu_i$. 
 The set
$\union_\nu \HIVE_{(2,1,0),(2,1,0)}^{\,\,\nu}$ is in
figure \ref{fig:minisamplehives} above.

Our goal is a self-contained proof of the following
positive formula for $\GLn$ tensor product multiplicities:

\begin{Theorem}\label{thm:main}
  Let $\lambda,\mu,\nu \in \integers^n_{\dec}$
  and let $V_\lambda,V_\mu,V_\nu$ be irreducible representations of $\GLn$
  with those high weights. Then the number of times $V_\nu$ appears
  as a constituent of the tensor product $V_\lambda\tensor V_\mu$
  is the number of lattice points in $\HIVE_{\lambda\mu}^\nu$.
\end{Theorem}

For example, figure \ref{fig:minisamplehives} is computing the tensor square 
$$ V_{(2,1,0)}^{\bigotimes 2} 
\iso V_{(4,2,0)}\oplus V_{(4,1,1)}\oplus V_{(3,3,0)}\oplus 
V_{(3,2,1)}^{\bigoplus 2} \oplus V_{(2,2,2)}.$$

While it doesn't make any sense to {\em count}
real-valued hives with fixed boundary (which is why we
insist on integer values), one can still consider the convex polytope
thereof, and relate it to the geometry of certain moduli spaces (see
the appendix to \cite{KTW}).  It is rather harder to formulate a ``real
version'' of skew Young tableaux!

\section{Recognizing the representation ring $Rep(\GLn)$}

\newcommand\Alt{\Lambda}

Recall that the representation ring $Rep(\GLn)$ has a basis
$\{[V_\lambda]\}$, $\lambda \in \integers^n_{\dec}$.
Let $\omega_i^n$ denote the ``fundamental weight'' 
$(1,\ldots,1,0,\ldots,0)$ with $i$ $1$s and $n-i$ $0$s, the
high weight of $\Alt^i \complexes^n$. (The notation is a little
nonstandard -- people usually just use $\omega_i$ -- but that would be
clumsy in lemma \ref{lem:strip} to come.)

The only other facts we will need about $Rep(\GLn)$ -- for which our
reference is \cite{FH} -- are that
\begin{itemize}
\item it is associative
\item it is generated by the fundamental representations 
  $[V_{\omega_i^n}]$ and $[V_{(-1,-1,\ldots,-1)}]$
\item $[V_\lambda] \, [\Alt^n \complexes^n]^{-1} = [V_{\lambda-(1,\ldots,1)}]$
  (we'll call this the {\dfn $\det^{-1}$ rule})
\item it satisfies the {\dfn Pieri rule}:
  $$ [V_\lambda] \, [\Alt^i \complexes^n] = 
       \Oplus_{\pi \in \{0,1\}^n,\, \sum\!\! \pi = i \atop
         \lambda + \pi \,\in\, \integers^n_{\dec}}
       V_{\lambda + \pi} $$ 
The sum is over those $0,1$-vectors $\pi$ with $i$ ones
(or equivalently those weights occurring in $\Alt^i \complexes^n$),
such that $\lambda+\pi$ is weakly decreasing.
\end{itemize} 

If $R$ is a ring-with-basis isomorphic to $Rep(\GLn)$, then it satisfies
the $\det^{-1}$ and Pieri rules; we now show that the converse is true.
(Essentially the same observation is used in \cite{T} 
and is surely much older.)

\begin{Proposition}\label{prop:ringiso}
  Let $R$ be a ring with $\integers$-basis 
  $\{b_\lambda\}, \lambda \in \integers^n_{\dec}$,
  satisfying the $\det^{-1}$ and Pieri rules. 
  Then the evident linear isomorphism 
  $\phi: Rep(\GLn) \to R$,  $[V_\lambda] \mapsto b_\lambda$ is also a
  ring isomorphism.
\end{Proposition}

\begin{proof}
  We want to show that $\phi(x y) = \phi(x)\phi(y)$. By linearity,
  it's enough to show it for $x$ a basis element $[V_\lambda]$.
  
  The Pieri and $\det^{-1}$ rules being true in both rings then tells
  us that this equation does hold if $x$ is a fundamental
  representation, $[V_{(1,\ldots,1,0,\ldots,0)}]$ or
  $[V_{(-1,-1,\ldots,-1)}]$. 

  More generally, let $y = b_{\mu_1} b_{\mu_2} \ldots b_{\mu_l}$ be
  a product of $l>0$ generators. Then
  $$ \phi(b_\lambda (b_{\mu_1} b_{\mu_2} \ldots b_{\mu_l}))
  = \phi((b_\lambda b_{\mu_1} b_{\mu_2} \ldots b_{\mu_{l-1}}) b_{\mu_l})
  = \phi(b_\lambda b_{\mu_1} b_{\mu_2} \ldots b_{\mu_{l-1}}) \phi(b_{\mu_l})
  $$
  and induction on $l$ takes care of the rest.
  (Note that the identity, $[V_{(0,\ldots,0)}]$, is itself a product
  $[V_{(1,\ldots,1)}] [V_{(-1,\ldots,-1)}]$ of two of our generators,
  so requiring $l>0$ does not cause us to miss this basis element.)

  So far we know that $\phi$ is establishing a ring isomorphism between the
  subspace of $Rep(\GLn)$ generated by the fundamental representations,
  and the image of that under $\phi$. But since the
  fundamental representations generate $Rep(\GLn)$,
  and $\phi$ is a linear isomorphism, that's actually a ring isomorphism
  between the two rings.
\end{proof}

In the rest of the paper our ring $R$ will be the {\dfn hive ring},
where the multiplication is {\em defined} by
$$ b_\lambda b_\mu = \sum_\nu \# \HIVE_{\lambda\mu}^\nu b_\nu. $$
The hardest part in applying proposition \ref{prop:ringiso}
will be to prove that $R$ is associative. Since we haven't proved
that yet it's a bit disingenuous to call it a ring, but we'll
do it anyway rather than having to rename it afterward.

Once we've checked $\det^{-1}$, Pieri, and associativity
for the hive ring, theorem \ref{thm:main} will follow from 
proposition \ref{prop:ringiso}.

\section{The hive ring satisfies the $\det^{-1}$ and Pieri rules}

\begin{Lemma}\label{lem:parallelograms}
  Let $p$ be a lattice parallelogram in the hive triangle $tri_n$, 
  with edges parallel to the edges in the triangular lattice, and
  $h$ a hive of size $n$. Then the sum of $h$'s entries at the two obtuse
  angles of $p$ is greater than or equal to 
  the sum of $h$'s entries at the two acute angles of $p$.
\end{Lemma}

\begin{proof}
  Add up all the rhombus inequalities from the rhombi inside and aligned with 
  $p$; everything cancels except the contributions from the four corners.
\end{proof}

\begin{Proposition}\label{prop:detinverse}
  In the hive ring, 
  $b_\lambda b_{(-1,\ldots,-1)} = b_{\lambda + (-1,\ldots,-1)}$.
  That is to say, the hive ring obeys the $\det^{-1}$ rule.
\end{Proposition}

\begin{proof}
  We're studying the hives with differences 
  $\lambda_i = h_{n-i,0,i} - h_{n-i+1,0,i-1}$ on the Northwest side, 
  and that are linear with slope $-1$ on the Northeast side
  (so $h_{0,n-z,z} = h_{0,n,0}-z$). We want to show there's exactly one,
  and it has $h_{i,0,n-i} = h_{i,n-i,0} - i$.

  Let $h \in \HIVE_{\lambda, (-1,\ldots,-1)}^\nu$ for some $\nu$.
  Consider the entry $h_{x,y,z}$, and the following two parallelograms
  in $tri_n$ with $[x,y,z]$ as a vertex:

  \centerline{\epsfig{file=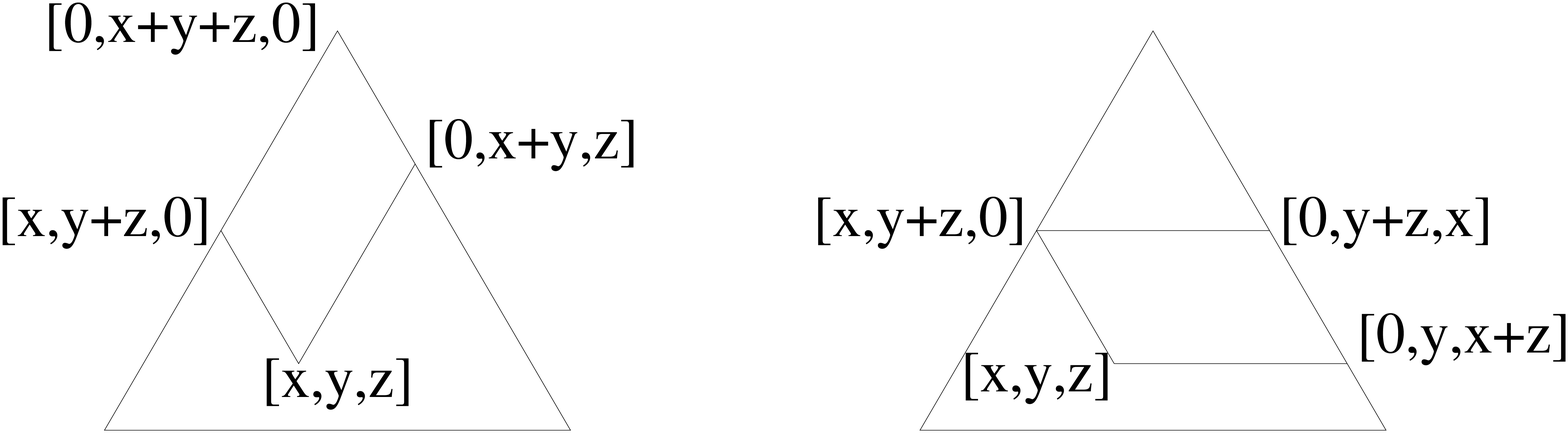,height=1.3in}}

  Let $\Lambda = \sum_i \lambda_i$ denote the value $h_{0,n,0}$ at the 
  top. Then the parallelogram inequalities of lemma \ref{lem:parallelograms},
  $$ h_{x,y+z,0} + h_{0,x+y,z} \geq h_{x,y,z} + h_{0,x+y+z,0} 
  \qquad \hbox{and} \qquad
  h_{x,y,z} + h_{0,y+z,x} \geq h_{x,y+z,0} + h_{0,y,x+z}, $$
  can be rewritten as
  $$ h_{x,y+z,0} + \Lambda - z \geq h_{x,y,z} + \Lambda 
  \qquad \hbox{and} \qquad
   h_{x,y,z} + \Lambda - x \geq h_{x,y+z,0} + \Lambda - x - z. $$
  These bound $h_{x,y,z}$ above and below by $h_{x,y+z,0} - z$. 

  In particular the South edge is given by $h_{x,0,z} = h_{x,z,0} - z$; 
  the only possible $h$
  has the differences $(\lambda_1 - 1, \lambda_2 - 1, \ldots, \lambda_n - 1)$
  across the bottom and the rest of the hive is uniquely determined.

  That shows uniqueness of the hive; how about existence? 
  The convexity of the function $h_{x,y,z} = h_{x,y+z,0} - z$ can
  be traced, with a bit of algebra, to the assumption that $\lambda$
  was weakly decreasing.
\end{proof}

This proposition can instead be proved by noting
that adding an linear function of $y$ and $z$ to a hive
produces a new hive, and by using the same inequalities to show that
$b_\lambda b_{\vec 0} = b_\lambda$, whose unique hive is constant on
NW/SE lines.

\begin{Lemma}\label{lem:strip}
  Let $h$ be an $n$-hive such that the differences down the NE edge are
  $\omega_i^n$. Then the differences down the strip one step in from
  the NE edge are either $\omega_i^{n-1}$ or $\omega_{i-1}^{n-1}$.
  
  Depending on which, the last difference $h_{1,0,n-1}-h_{0,0,n}$
  across the bottom either agrees with the last difference
  $h_{1,n-1,0} - h_{0,n,0}$ on the NW side, or is one larger, respectively.
\end{Lemma}

\begin{proof}
  For short, write $h_{1,n-1,0}$, $h_{0,n,0}$, $h_{1,0,n-1}$, $h_{0,0,n}$
  as $x$, $x+a$, $y$, $x+a+i$ respectively. (That $h_{0,0,n}=x+a+i$ follows
  from the assumption that the differences across the NE side are
  $\omega_i$, which has total $i$.)

  \centerline{\epsfig{file=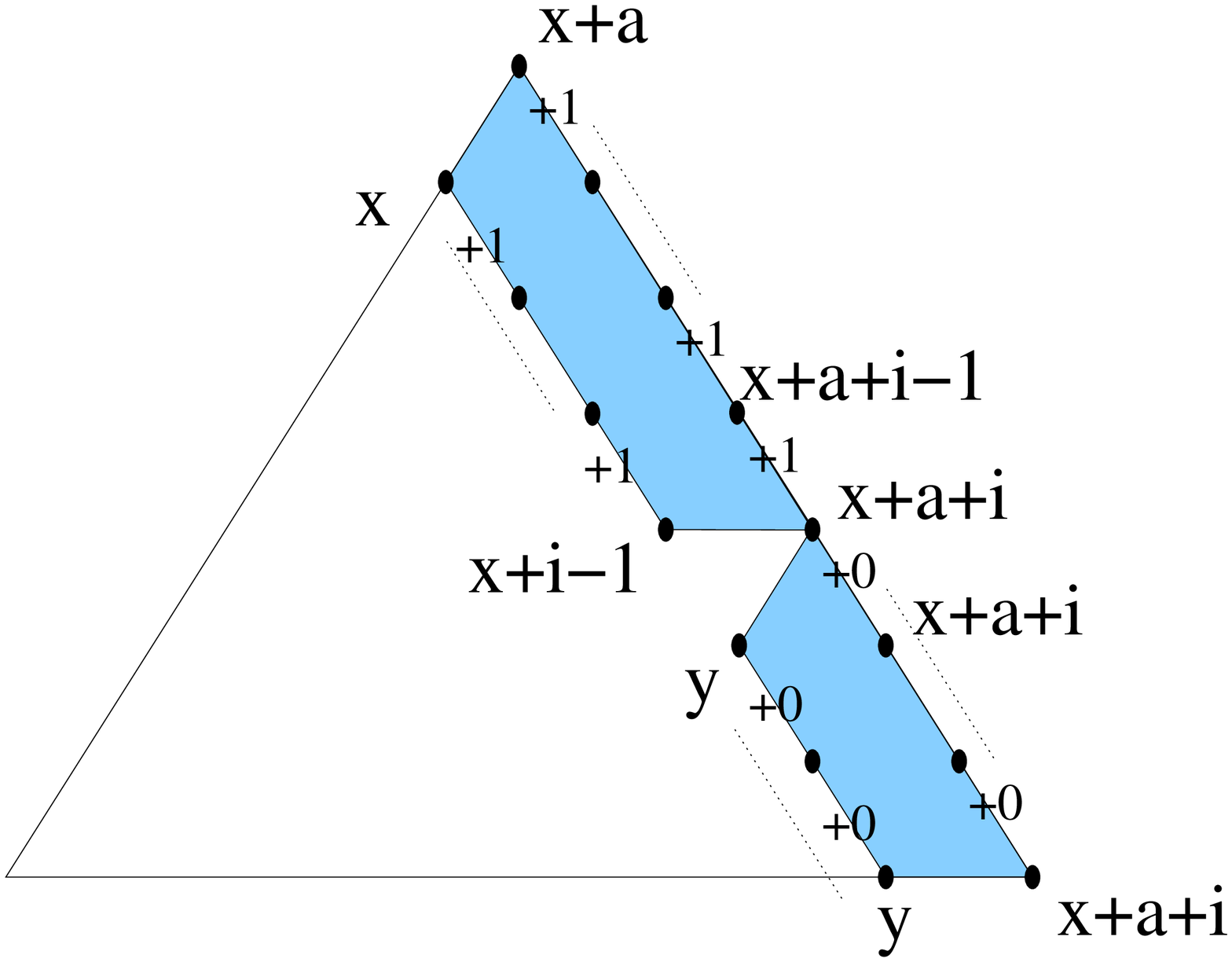, height=2in}}
  
  Using only the rhombus inequalities in the shaded regions of the
  figure above, and the same line of argument as in proposition
  \ref{prop:detinverse}, we can show that $h_{1,n-1-i,i} = x+i$ and
  $h_{1,n-i-2,i+1} = y$. (These are the two adjacent interior entries
  indicated in the figure.)

  Now the two rhombus inequalities relating those hive entries and the
  NE boundary say $ x+i \geq y \geq x+i-1 $. In particular, the 
  difference $(x+a+i) - y$ is either $a$ or $a+1$, and that binary choice
  determines the rest of the strip.
\end{proof}

\begin{Proposition}\label{prop:pieri}
  In the hive ring, 
  $$ b_\lambda b_{\omega_i} 
  = \sum_{\pi \in \{0,1\}^n,\, \sum\!\!\pi = i \atop 
    \lambda + \pi\, \in\, \integers^n_{\dec}}
  b_{\lambda + \pi}. $$
  In other words, ``the hive ring obeys the Pieri rule''.
\end{Proposition}

\begin{proof}
  Let $h$ be a hive with differences $\lambda$ on the NW side, 
  $\omega_i$ on the NE side. Rip off the NE strip from it and repeat,
  each time producing a hive one size smaller. 

  By inductive use of lemma \ref{lem:strip}, we see that the differences
  on the NE side go from $\omega_i$ (at size $n$) to $\omega_0$ (at size $0$),
  so the differences across the bottom agree with $\lambda$ in $n-i$ places
  and are one larger in $i$ places. Moreover, the hive is uniquely
  determined by its labels on the bottom edge.

  By proposition \ref{prop:dec}, the differences in the labels 
  across the bottom are still decreasing. 
  This, plus the previous paragraph, establishes the Pieri rule
  as an upper bound.

  Given a $0,1$-string $\pi$ with $i$ ones such that $\lambda+\pi$
  is dominant (and so {\em should} be giving a term in the Pieri rule),
  we can glue together the strips from lemma \ref{lem:strip} and hope
  that we get a hive. The only rhombus inequalities left to check are
  those intersecting two adjacent strips, and we leave this to the reader.  
\end{proof}

\section{The hive ring is associative}\label{sec:assoc}

First off, what's the equation we're trying to prove? 
Let $h_{\lambda\mu}^\sigma = \#\HIVE_{\lambda\mu}^\sigma$,
the structure constant in the hive ring. Then
$$ (b_\lambda b_\mu) b_\nu 
= \sum_\sigma h_{\lambda\mu}^\sigma b_\sigma b_\nu
=\sum_\sigma \sum_\pi h_{\lambda\mu}^\sigma h_{\sigma\nu}^\pi b_\pi
$$
whereas
$$ b_\lambda (b_\mu b_\nu) = \sum_\tau  b_\lambda h_{\mu\nu}^\tau b_\tau
= \sum_\tau \sum_\pi h_{\mu\nu}^\tau h_{\lambda\tau}^\pi b_\pi $$
Comparing coefficents of $b_\pi$, we see that we need to prove
$$  (*)\qquad\qquad\qquad\qquad\qquad\qquad
\sum_\sigma  h_{\lambda\mu}^\sigma h_{\sigma\nu}^\pi 
= \sum_\tau  h_{\mu\nu}^\tau h_{\lambda\tau}^\pi. 
\qquad\qquad\qquad\qquad\qquad\qquad$$

Consider a tetrahedron balanced perfectly on an edge, from directly above;
the boundary of what you see is a square. Label the edges of this square 
(starting from the top left vertex and going clockwise)
with the partial sums of $\lambda,\mu,\nu,\pi^*$.
(The dominant weight $\pi^*$ is $(-\pi_n, -\pi_{n-1},\ldots,-\pi_1)$,
the highest weight of the contragredient representation $(V_\pi)^*$.
One could say it comes up because we're reading that edge of the
hive backwards.)

If the top edge is labeled $\sigma$, then the number of ways
of labeling the upper two faces with hives is
$h_{\lambda\mu}^\sigma h_{\sigma\nu}^\pi$. 
Without fixing the labeling on that top edge,
it's $ \sum_\sigma  h_{\lambda\mu}^\sigma h_{\sigma\nu}^\pi$.
The corresponding statement for the lower two faces gives the other sum.

\begin{center}
  \epsfig{file=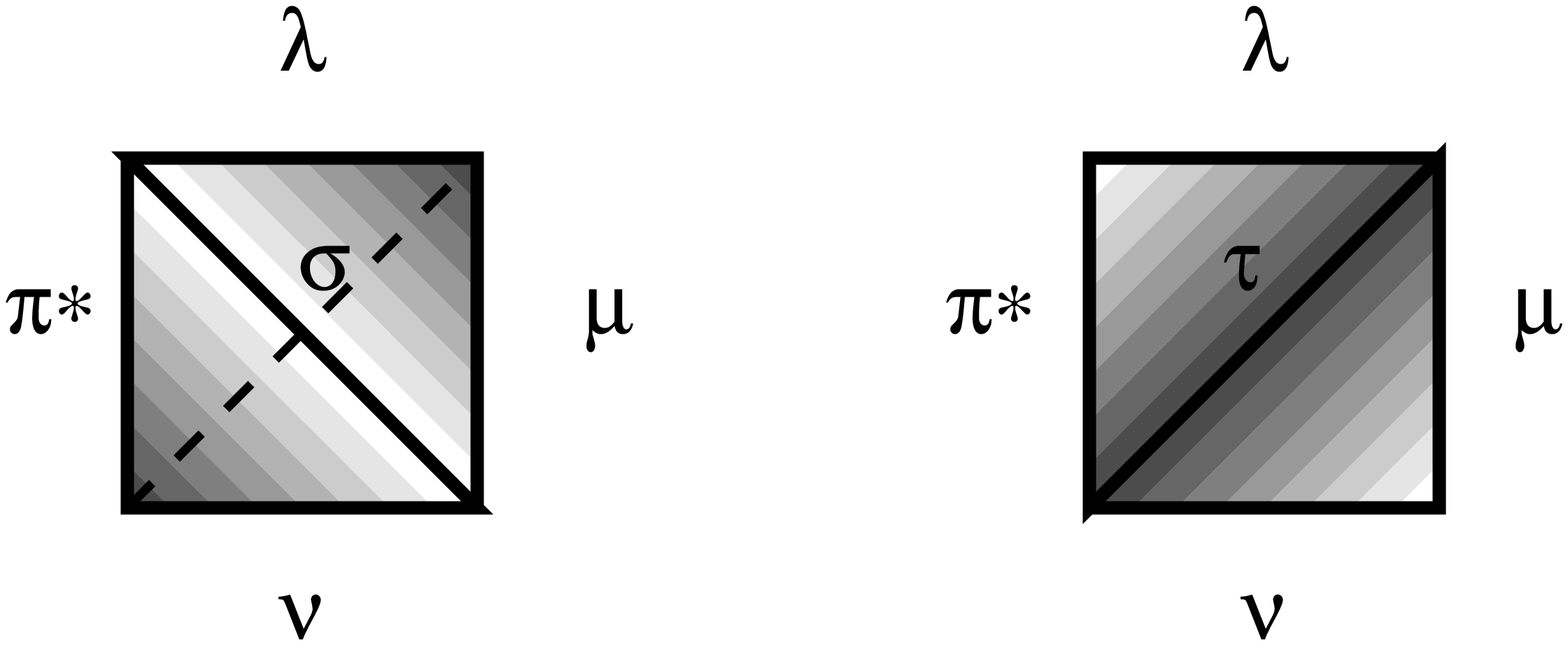, height=1.5in}
\end{center}

\begin{Theorem}\label{thm:excavation}
  There is a continuous, piecewise linear bijection between ways of
  labeling the upper two faces of this tetrahedron with a pair of 
  {\em real} hives and ways of labeling the lower two faces, with
  given fixed labels $\lambda,\mu,\nu,\pi^*$ around the four
  non-horizontal edges.

  Moreover, each formula for a label on a bottom face is a
  ``tropical Laurent polynomial'' in the entries on the top two faces,
  meaning it can be written as a maximum over some linear forms.

  This bijection on matched pairs of real hives restricts to a bijection
  on matched pairs of integral hives, which establishes equation $(*)$ above.
\end{Theorem}

\begin{proof}
  This tetrahedron of size $n$ breaks up into little tetrahedra, little
  upside-down tetrahedra, and octahedra (think about the $n=2$ case).
  In coordinates, let $tet_n = \{[x,y,z,w]\in \naturals^4 : x+y+z+w = n\}$.
  Then the right-side-up tetrahedra have vertices
  $$ [x+1,y,z,w], [x,y+1,z,w], [x,y,z+1,w], [x,y,z,w+1], $$
  the octahedra have vertices
  $$ [x\!+\!1,y\!+\!1,z,w],[x\!+\!1,y,z\!+\!1,w], [x\!+\!1,y,z,w\!+\!1], 
        [x,y\!+\!1,z\!+\!1,w], [x,y\!+\!1,z,w\!+\!1], [x,y,z\!+\!1,w\!+\!1], $$
  and the upside-down tetrahedra have vertices
  $$ [x+1,y+1,z+1,w], [x+1,y,z+1,w+1], [x+1,y+1,z,w+1], [x+1,y+1,z+1,w]. $$

  Imagine the tetrahedron as initially being ``full'' of these pieces,
  which we will remove one by one from above,
  each being removable only when everything above is already out of the way.
  Along the way, we'll label all the interior lattice points with
  numbers. When we're done, leaving only the bottom two faces, it will
  turn out that we have two hives there.

  Whenever we remove a little tetrahedron, we don't expose any new lattice
  points. Whenever we remove an octahedron, though, one of the old vertices
  (a local height max) goes with it and a new one becomes visible (a local 
  height min).
  As we go, we label the vertices exposed according to the following formula:
  $$ e' := \max(a+c,b+d) - e $$
  where $e$ was the label at the top, and $a,b,c,d$ the labels around
  the equatorial square. Our references for this {\dfn octahedron rule} 
  are \cite{P,FZ} (though it is much older, such as in \cite{RR}).

  When we're done, we have labeled the bottom two faces. The process ---
  which we call the {\dfn excavation} of $tet_n$ ---
  obviously provides its own inverse (the equation above is symmetric
  in $e$ and $e'$), and preserves integrality.

  It remains to see that what we get on the bottom is a pair of hives,
  i.e. satisfies the rhombus inequalities. We will show now that {\em every}
  unit rhombus in the tetrahedron gives a true rhombus inequality.

  Say we've partially excavated, and every rhombus above the level
  so far dug out has satisfied this inequality. Now we extract a piece;
  this exposes some new rhombi that we need to check. 

  {\em The $n=2$ case.}
  We remove the top two tetrahedra, then the octahedron, then a bottom
  tetrahedron. From the top, we see the labels\\
  \centerline{\epsfig{file=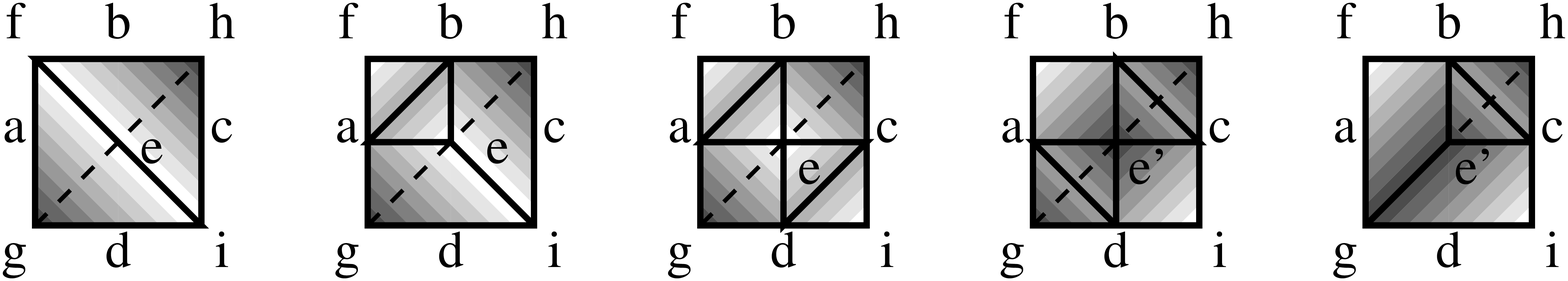, width=6in}} \\
  where the heavy (resp. dotted) lines indicate visible (resp. hidden) creases,
  and the shading indicates depth. From the South-Southeast ($d$ in
  front, $b$ in back), the process looks like this:\\
  \centerline{\epsfig{file=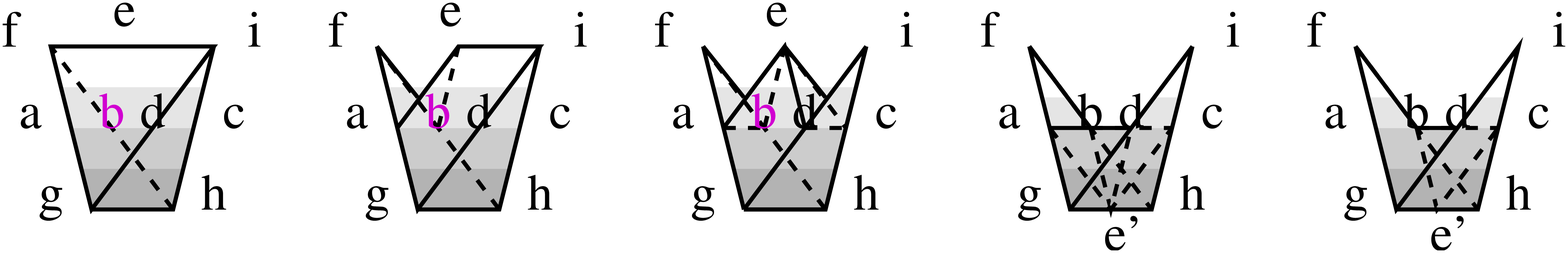, width=6in}} 
  and at the end only the $bce'h$ tetrahedron is left.

  The first two moves, removing the $abef$ and $edci$ tetrahedra, 
  expose no new lattice points (only the creases change).
  The next move exposes the $e'$ lattice point, 
  and thus the rhombus with obtuse vertices $a,b$,
  acute $f, e' = \max(a+c,b+d)-e$. We want to show that
  $$ a+b \geq f+\max(a+c,b+d)-e $$
  or equivalently
  $$ a+b \geq f+a+c-e,\qquad a+b \geq f+b+d-e $$
  which follow from the $b+e\geq f+c$, $a+e\geq d+f$ inequalities on the top.

  The third move exposes the rhombus with obtuse vertices $a,e'$,
  acute $b,g$. We want to show that
  $$ a+\max(a+c,b+d)-e \geq b+g $$
  so it's enough to show one of them: $a+b+d-e\geq b+g $. 
  This follows from $a+d\geq e+g$ on the top.

  While we haven't explicitly handled all the rhombi in this size
  $2$ tetrahedron -- or even finished excavating it; the $bhce'$ tetrahedron
  is still in place -- the other rhombi are equivalent to these two under 
  the evident $Z_2\times Z_2$ symmetries.

  {\em The general case.} Any rhombus exposed fits into a size $2$ 
  tetrahedron, so we just have to apply the $n=2$ case over and over.

  Finally, we need the ``tropical Laurent polynomial'' statement.
  The rule $e' = \max(a+c,b+d)-e$ is the {\em tropicalization} of the
  subtraction-free
  rational function $E' = (AC + BD)/E$, meaning that $+,\times,/\,$
  have been replaced with $\max,+,-$. As a very special case of 
  the main theorem 1.6 in \cite{FZ}, one knows
  that if one uses this rational function recurrence during excavation,
  the labels on the bottom two faces are {\em Laurent polynomials} in
  the labels on the top two (rather than merely rational functions as
  one would expect). 
\end{proof}

Feeding this theorem and propositions \ref{prop:detinverse} 
and \ref{prop:pieri} into proposition \ref{prop:ringiso}, 
we obtain theorem \ref{thm:main}.

Since this paper was first written, it was proven in 
\cite{S} that the coefficients in these Laurent polynomials are all $1$,
and the monomials identified.
We state this result, in the special case relevant here,
in section \ref{sec:speyer}.

\section{The honeycomb interpretation: scattering}\label{sec:honey}

This section is distinctly less detailed than the others, and is largely
for motivation. We recall briefly the {\em honeycombs} of \cite{KT1}, 
which are in $1:1$ correspondence with hives, but better suited for
some aspects of their study.

Let $\B$ denote the plane of triples of real numbers with zero sum.
Define the {\dfn coordinate directions} in $\B$ to be parallel to
$(0,1,-1)$, which we will draw as Northwest, $(-1,0,1)$, which we
will draw as Northeast, and $(1,-1,0)$, which we will draw as South.
A line segment oriented parallel to a coordinate direction has a
{\dfn constant coordinate}, the one of the three coordinates
constant along the edge.

\begin{figure}[htbp]
  \begin{center}
    \leavevmode
    \epsfig{file=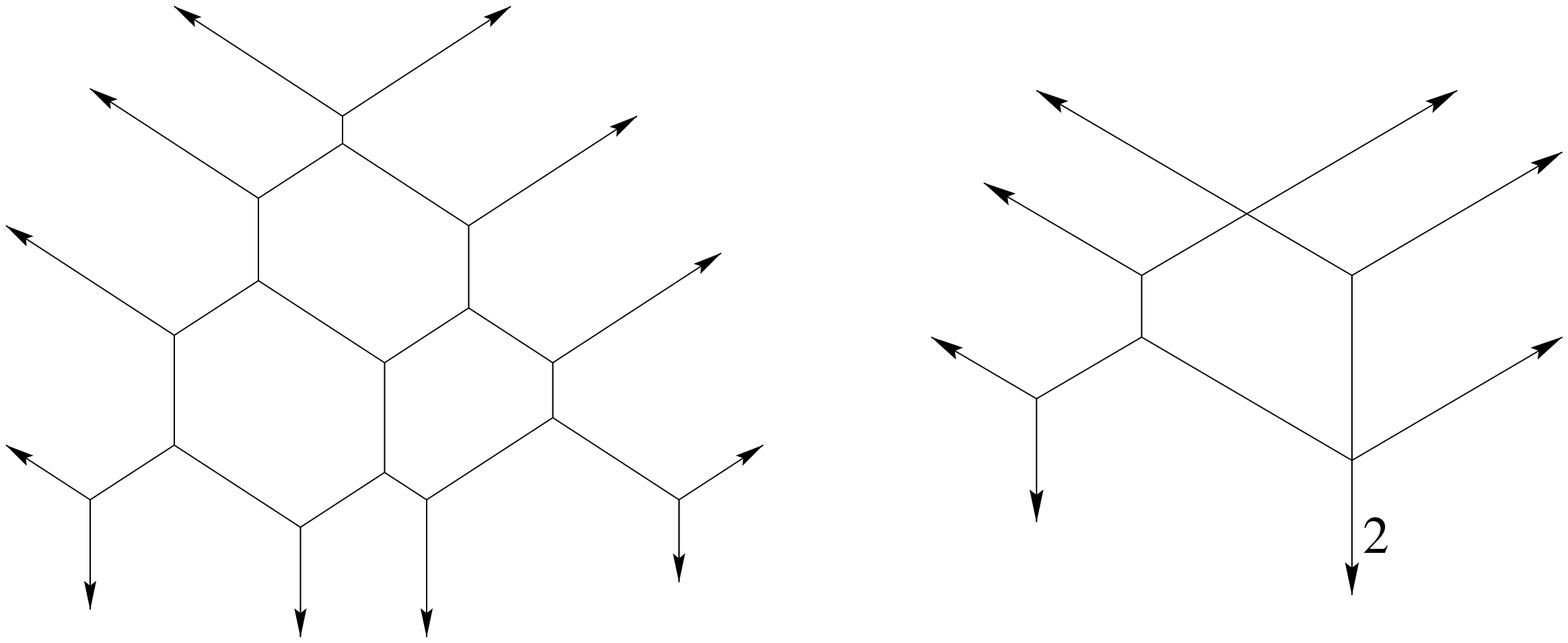,height=2in}
    \caption{Two honeycombs, of sizes $4$ and $3$.
      The left one is more typical in having only {\bf Y} vertices.
      All edges are multiplicity $1$, except for the edge labeled $2$ 
      in the right-hand honeycomb.}
    \label{fig:honeyex}
  \end{center}
\end{figure}

A {\dfn honeycomb} is a measure on $\B$, constructed by summing the
Lebesgue measure on a finite number
of coordinate-oriented line segments (which may be unbounded), such that
\begin{itemize}
\item each unbounded ray goes in a coordinate direction (not its negative)
\item around each point, the total ``pull'' of the up-to-six edges
  emanating from that point is the zero vector.
\end{itemize}
Note that some of the segments may overlap (or even coincide),
leading to multiplicities along the edges. Two honeycombs
are displayed in figure \ref{fig:honeyex}.

Since any vertex of a honeycomb satisfies this ``zero-tension'' condition,
the honeycomb as a whole does so (by a sort of Green's theorem), 
so the number of edges (counted with multiplicity) emanating in the Northwest, 
Northeast, or South directions must be the same number $n$. We call this
the {\dfn size} of the honeycomb.

From a hive of size $n$, we construct a honeycomb of size $n$ as follows.
There is one honeycomb edge for each unit edge connecting two vertices
in $tri_n$, but perpendicular to it (living in the dual graph).
The constant coordinate on that honeycomb edge is the difference
of the two labels in the hive, up to a certain sign. To determine
this sign, look for the unit triangle $\Delta$ in $tri_n$ aligned with
$tri_n$, and containing the two hive labels and an extra vertex.
The constant coordinate assigned is then the label on $\Delta$
counterclockwise of the extra vertex, minus the label on $\Delta$
clockwise of the extra vertex.

The vertices of the honeycomb then correspond to the linear regions in
the hive. The rhombus inequalities on the hive, reinterpreted,
state that the edges of
the honeycomb are of nonnegative length. It is quite tricky to prove that
this map from hives to honeycombs is in fact a bijection (theorem 1
of \cite{KT1}).

We are now in a position to describe ``honeycomb scattering'', 
a honeycomb interpretation of the tetrahedron-evacuation bijection 
from section \ref{sec:assoc}. This was the form in which one of us (CW)
first found this proof of associativity.

\begin{figure}[htbp]
  \centering
  \epsfig{file=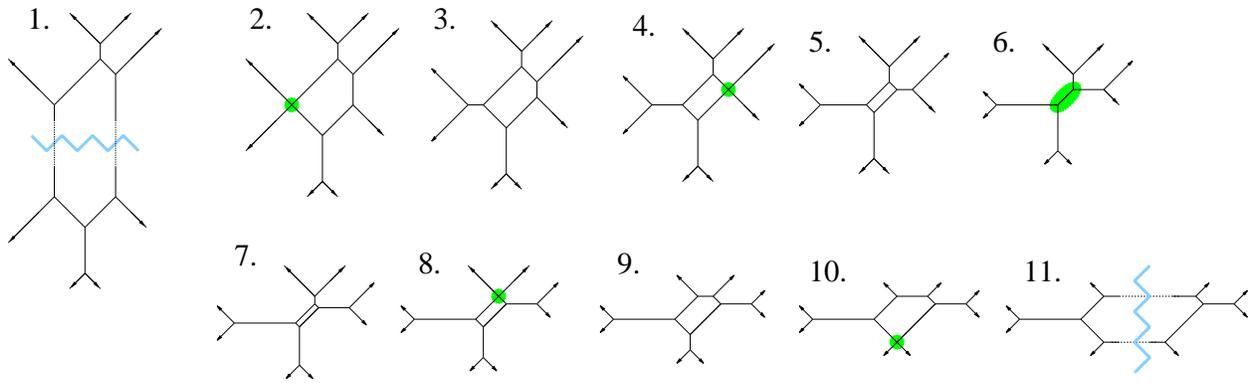, width=6.5in}
  \caption{The eleven stages of two size $2$ honeycombs scattering off 
    one another. The collisions are circled. Before the double
    collision in (6), the rectangle is shrinking; after that collision, 
    the rectangle bounces back out.}
  \label{fig:scatterex}
\end{figure}

Let $\HONEY_{\lambda\mu}^\nu$ denote the set of honeycombs whose 
boundary edges have constant coordinates $\lambda$ in the Northwest direction,
$\mu$ in the Northeast direction, and $\nu$ in the South direction.
To prove 
$$ \sum_\sigma  h_{\lambda\mu}^\sigma h_{\sigma\nu}^\pi 
= \sum_\tau  h_{\mu\nu}^\tau h_{\lambda\tau}^\pi, $$ 
consider the set of pairs of honeycombs 
$$ (h,h') \in \union_\sigma 
\big( \HONEY_{\lambda\mu}^\sigma \times \HONEY_{\sigma\nu}^\pi \big). $$
We can draw such a pair $(h,h')$ by rotating $h'$ by $180^\circ$, 
translating it some large distance in the $(1,-1,0)$ direction, and
gluing it to $h$, as in the first entry in
figure \ref{fig:scatterex}. For reasons to be
explained below we also compress this picture in the $y$ direction,
making the edges all $90^\circ$ and $45^\circ$ from one another
rather than $60^\circ$.

Now pull the lower honeycomb $h'$ upwards, while sending the upper
one downwards, at the same constant speed
(this interpretation involves a time coordinate).
At some point a vertex of $h'$ will collide with one of $h$.
We give an example of the whole process in figure \ref{fig:scatterex},
and (before they are given below)
we invite the reader to guess the general rules defining scattering.

The rules for scattering are as follows. 
Each vertex in $h$ (respectively, $h'$) is
given an initial velocity of down (respectively, up).
All vertices move until there is a collision. 

Generically, the first sort of collision met is
that of a single edge contracting (as in (1-2), 
(3-4), (7-8), (10-11) of figure \ref{fig:scatterex}). When this happens,
we redirect the vertically colliding particles to move left and right,
conserving total energy and momentum. Note that what used to be a vertical
line connecting them is now horizontal.

After some of these collisions, a second type of collision is possible,
involving two pairs of converging vertices making a rectangle collapse
(as in (5-6) of figure \ref{fig:scatterex}). These are again redirected out,
conserving energy and momentum.

(For generic $h$ and $h'$, these two are the only sorts of collisions
that can happen during the whole scattering. Nongeneric $h,h'$ can be
understood by taking limits from the generic case, so we won't dwell
on them.)

It is true at the beginning, and remains true after either type of
collision, that if a vertex is attached to a vertical line then it
is moving vertically and will eventually get into a 2-vertex collision.
Note that each 2-vertex collision increases the number of particles
moving horizontally and decreases the number of vertical edges, 
and each 4-vertex collision decreases the number of
particles moving directly toward one another. 

So there are only a finite 
number of collisions; when the scattering is over, all particles are moving
horizontally, there are no vertical edges, and all the left-moving
particles are left of all the right-moving particles. At this point
we can cut the diagram in half along the growing edges, and we get
two new squashed honeycombs -- except that they've been squashed along
the $x$ direction rather than the $y$.

\subsection{Honeycomb scattering vs. hive excavation.}

These two types of collisions -- 2-vertex and 4-vertex -- correspond to
the tetrahedral and octahedral excavations in theorem \ref{thm:excavation}.
The hive labels are linearly related to the constant coordinates on the
diagonal lines. Note that during a 2-vertex collision, the four diagonal
lines incident move at a constant speed -- on the excavation side, this
reflects the fact that excavating a tetrahedron exposes no new vertices
and requires no new labeling. The $\max$ of the octahedron rule is
implemented by the two ways a rectangle can collapse -- whichever one
comes first determines how the vertices bounce back out.

In the hive excavation picture, it is easier to deal with degenerate cases
uniformly -- the octahedron rule still applies, it just happens that the
$\max$ involved is achieved twice. Also the hive picture doesn't introduce
this spurious ``time'' coordinate; in particular, the excavation of the
large tetrahedron can be done in many different ways, all giving the
same answer.

On the other hand, in the honeycomb picture it is more manifest that the
limiting object after scattering is again two honeycombs glued together,
rather than having to check the rhombus inequalities in various cases,
as occurred in theorem \ref{thm:excavation}.

\subsection{The scattering rule in \cite{GP}.}

A very similar result is proven in \cite{GP}, though phrased in terms
of braid relations rather than excavation. It is less immediately
obvious that their rule is constructing an associator, since it 
involves {\em six} inputs rather than four. In fact their pseudo-line
arrangements can be corresponded to a partially excavated {\em cube}
of size $n$, rather than our tetrahedron, and their braid move
corresponds to the removal of one little cube. 


\section{A closed form for the associativity bijection}\label{sec:speyer}

As promised in section \ref{sec:assoc}, we give a closed-form expression
for each entry of the bottom two hives as a ``tropical Laurent polynomial''
(a single maximum over a family of linear expressions) in the entries
of the top two. 

This is a special case of a general result proved by Speyer in \cite{S}, 
conjectured in
\cite{FZ}: during any order in which we excavate a tetrahedron, a label
$b$ exposed at time $t_2 > t_1$ is a Laurent polynomial 
with positive coefficients 
in those labels on the surface at time $t_1$. 
Speyer proves a more precise conjecture due to Propp:
each Laurent monomial in $b$ corresponds to a {\dfn matching} 
of a certain graph $G_b$,
meaning a subset of the edges covering each vertex exactly once,
where the graph $G_b$ is determined by the $t_1$-surface and the $t_2$-entry.
Intriguingly, the ``graphs with open faces'' $G_b$ that Speyer
constructs to prove this look like partially-scattered honeycombs!

\begin{figure}[htbp]
  \centering
  \epsfig{file=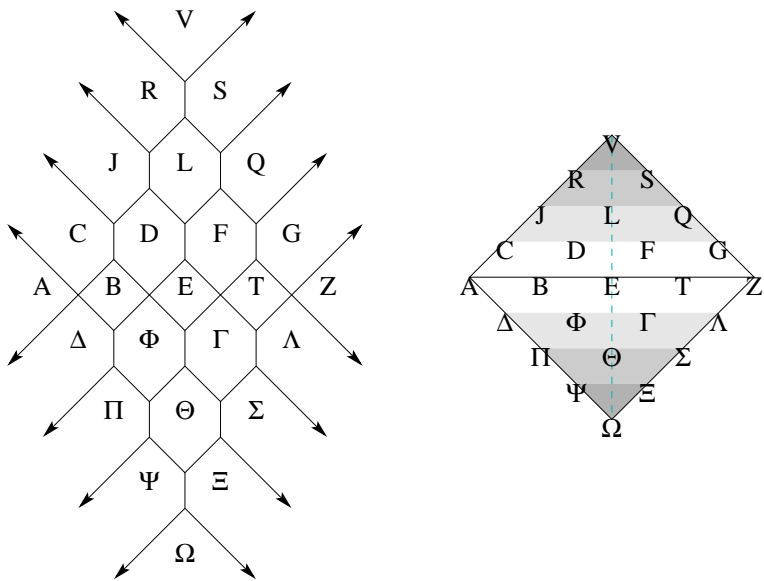,height=3in}
  \caption{Two standard $4$-honeycombs caught at the moment of 
    first scattering (left). The regions correspond to labels on the top
    two faces of a tetrahedron (right), as indicated.}
  \label{fig:speyer1}
\end{figure}

Consider the graph $G$ (with some unbounded edges) constructed by scattering
two {\em standard} $n$-honeycombs (meaning, all finite edges of length $1$) off
one another, stopping exactly when the first $n$ collisions occur 
(simultaneously). This graph has $n-1$ rhombi and two triangular arrays 
of ${n-1 \choose 2}$ hexagons. The $n=4$ example is drawn in figure 
\ref{fig:speyer1}. Its regions correspond to the labels on the
top two faces of a tetrahedron to be excavated.

To expose a bottom entry $b$, there is a unique minimal set of unit
tetrahedra and octahedra that must be excavated. For example, to expose
the bottom entry below the $F$ in figure \ref{fig:speyer1}, one must
remove the octahedra with top points labeled by $T,E,F$. 

\begin{figure}[htbp]
  \centering
  \epsfig{file=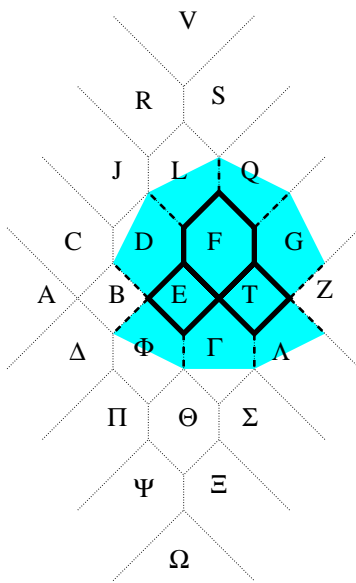,height=3in}
  \caption{The graph with open faces whose matchings will calculate
    the entry directly below the label $F$, with three interior and
    seven exterior faces.}
  \label{fig:bgraph}
\end{figure}

We are now ready to describe the graph $G_b$ that Speyer associates to a bottom
entry $b$: it is the minimal subgraph of $G$ enclosing those entries that 
must be excavated to get to $b$. In fact one must also consider the
adjoining hexagonal regions in Speyer's definition of $G_b$ as a ``graph
with open faces''. See figure \ref{fig:bgraph} for the
case $b$ being the entry below the $F$.

\begin{figure}[htbp]
  \centering
  \epsfig{file=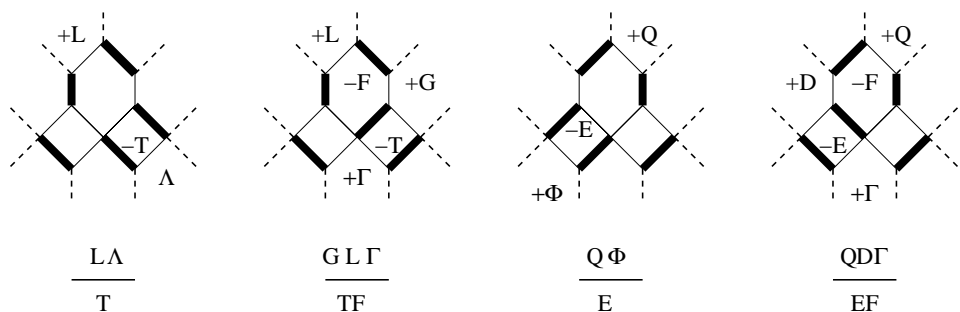,width=5in}
  \caption{The matchings of $G_b$ --- in each figure, each vertex touches
    exactly one heavy edge --- and their matching monomials.}
  \label{fig:bmatchings}
\end{figure}

Given such a graph $G_b$, and a matching $\mu$ of it, Speyer defines the 
{\dfn matching (Laurent) monomial} $m_\mu$ as the product over all faces
of $G_b$ (including the exterior hexagons), of the corresponding variable
raised to the power
\begin{itemize}
\item one minus the number of adjacent edges in $\mu$, for a rhombus
  or external hexagon, or
\item two minus the number of adjacent edges in $\mu$, for an
 interior hexagon.
\end{itemize}
In figure \ref{fig:bmatchings} we draw all the matchings of our
$G_b$ from figure \ref{fig:bgraph}, and compute the matching monomials.

Speyer's theorem (for our case) now reads as follows:

\begin{Theorem*}\cite{S}
  Let $b$ be a label on the bottom of the hive tetrahedron, and $G_b$
  the minimal subgraph of $G$ enclosing the entries on the top that
  must be excavated to expose $b$.

  Then if during excavation, we use the rational function octahedron
  recurrence $E' = E^{-1} (AC+BD)$, the resulting value of $b$ can be
  computed as a sum over all matchings $\{m\}$ of $G_b$, of the
  associated matching monomials $\{\mu_m\}$.

  If we instead use the tropical recurrence, the value of $b$ can be
  computed as a maximum over all matchings, of the corresponding
  linear forms.  
\end{Theorem*}

To compute the figure \ref{fig:bgraph}
example directly: the octahedron recurrence gives
\begin{align*}
  T' &= T^{-1} (F \Lambda + G \Gamma) \\ 
  E' &= E^{-1} (F \Phi + D \Gamma) \\
  F' &= F^{-1} (T' L + Q E') \\
     &= F^{-1} (T^{-1} (F\Lambda + G\Gamma) L + Q E^{-1} (F\Phi + D\Gamma))\\
     &= T^{-1} \Lambda L + F^{-1} T^{-1} G\Gamma L 
        + Q E^{-1} \Phi + F^{-1} Q E^{-1} D\Gamma
\end{align*}
which agrees with the theorem, being the sum of the terms from
figure \ref{fig:bmatchings}. The tropical version is therefore
$$
  F' = \max \{ -T + \Lambda + L, - F - T + G + \Gamma + L,
         Q - E + \Phi, - F + Q - E + D +\Gamma).
$$

\bibliographystyle{alpha}

\end{document}